\RequirePackage{fix-cm}
\documentclass[smallextended]{svjour3}       % onecolumn (second format)
\smartqed  % flush right qed marks, e.g. at end of proof
\usepackage{graphicx}
\usepackage{amssymb}
\begin{document}

%\title{Extended Dijkstra
%algorithm  and  Moore-Bellman-Ford algorithm
% \thanks{This work/[Author name] was supported by [institutions/academy/funds] (No. xxxxxx)}
%}
%\subtitle{Do you have a subtitle?\\ If so, write it here}

%\titlerunning{Short form of title}        % if too long for running head

%\author{$\textbf{Cong-Dian Cheng}^{1}$
% \and
%  $\textbf{Second Author}^{2}$ %etc.
%}

%\authorrunning{Short form of author list} % if too long for running head

%\institute{ F. Author \at
 %   \email{zhiyang918@163.com}
      %  \\
    %             \emph{Present address:} of F. Author  %  if needed
     %\and
    % S. Author \at
    % \email{fauthor@example.com} \\
%    \\
%1   first adress(School of Intelligence Technology, Geely University
%of China, Chengdu 641400, PR China)
%\\
%2   second adress
%}

%\date{Received: date / Accepted: date}

%\maketitle

\begin{center}\textbf{\large{Extended Dijkstra
algorithm  and  Moore-Bellman-Ford algorithm}} \\~\\
\small{Congdian Cheng\\~\\School of Intelligence Technology, Geely
University of China,
 Chengdu 641400, PR China\\\textbf{E-mail}: zhiyang918{\rm
 @}163.com.}

\end{center}

\begin{abstract}
Study a general single-source shortest path  problem, which is
motivated by current interest in needing to extend the total weight
function of paths on a network and the classical shortest path
problem. Firstly, define the path functional on a set of certain
paths with same source on a graph; introduce a few concepts
 of the defined path functional;
and make some discussions on the properties of the path functional.
Secondly, develop a kind of general single-source shortest path
problem (GSSSP). Thirdly, following respectively the approaches of
the well known Dijkstra's algorithm and Moore-Bellman-Ford
algorithm, design an extended Dijkstra's algorithm (EDA) and an
extended Moore-Bellman-Ford algorithm (EMBFA) to solve the problem
GSSSP under certain given conditions. Fourthly, under the assumption
that the value of related path functional for any path can be
obtained in $M(n)$ time, prove respectively the algorithm EDA
solving the problem GSSSP in $O(n^2)M(n)$ time and the algorithm
EMBFA solving the problem GSSSP in $O(mn)M(n)$ time. Finally, some
applications of the designed algorithms are shown with a few
examples. What we done can improve both the researches and the
applications of the shortest path theory.
% \PACS{PACS code1 \and PACS code2 \and more}

\end{abstract}
\section{Introduction}
\label{intro} %\uppercase\expandafter{\romannumeral1}
%\romannumeral2
Shortest path problems are the best known class of combinatorial
optimization problems and have been extensively studied for more
than half a century, which have many applications in
network,electrical routing, transportation, robot motion planning,
critical path computation in scheduling, quick response to urgent
relief, etc; and can also unify framework for many optimization
problems such as knapsack, sequence alignment in molecular biology,
inscribed polygon construction, and length-limited Huffman-coding,
etc. For the basic knowledge of shortest path problems, please refer
to  chapter 7 of the monograph \cite{biblabel[Korte]} (Korte and
Vygen 2000) and the other literatures afterword.

The classical single-source shortest path problem of network,
denoted by CSSSP, is the most famous   one of shortest path
problems, and a lot of works have been done to study and solve this
kind of shortest path problem. Among many algorithms for the problem
CSSSP, Dijkstra's Algorithm (DA) and Moore-Bellman-Ford Algorithm
(MBFA), called also by Bellman-Ford Algorithm, are two well-known
and most fundamental, which have now been the core technique to
solve many optimization problems. As we all know, the first one can
solve the  problem CSSSP with nonnegative edge weights in $O(n^2)$
time and the second one can deal with the problem CSSSP with
arbitrary conservative weights in $O(nm)$  time. Here $n$ and $m$
denote respectively the number of vertices and the number of edges
on the underlying graph. See, e.g., chapter 7 of the monograph [11]
(Korte and Vygen 2000) and the literatures [2,3,8,15] (Bellman 1958,
Dijkstra 1959, Ford 1956, Moore 1959).

The study of shortest paths is a research area with a long history.
And there are still substantial research works in the area from this
century. On the work over the last two decades,  please see the
following introduction for example.

Motivated by the recent interest in pricing networks and other
computational problems, Hershberger and Suri [9] (2001) studied the
problem how to determine the Vickrey payments of all the agents
between the given two nodes with a pricing network in less time,
which is closely related to the shortest path algorithm and the
replacement path problem, see Hershberger et al. [10] (2003); and
they proposed an algorithm to complete the computation in
essentially the same time bound as one single-source shortest path
computation. Hershberger et al. [10] (2003) also explained and
investigated the replacement path problem and the other shortest
path problems; and made some results on the time complexity of
computing the related problems. Du et al. [6] (2008), based on the
known heuristic algorithm CST for Euclidean Steiner Tree (EST),
proposed the algorithms CST(A) to find EST, and by making the
worst-case analysis of algorithms CST and CST(A), presented the
restricted submodularity technique to analyze approximation
algorithm with nonsubmodular functions; for Connected Dominating Set
(CDS), based on the known $1$ - greedy algorithm, they also proposed
$(2k-1)$ - greedy algorithm to compute a minimum CDS, and by making
the worst-case analysis of $(2k-1)$- greedy algorithm, presented the
shifted submodularity technique to analyze approximation algorithm
with nonsubmodular functions; in the process they obtain some
excellent results.

Noted the fact that  the transportation system of a city and the
roads of the city can be respectively modeled by a network and its
edges, some roads of which may be blocked at certain times; and the
traveler only observes that upon reaching an adjacent site of the
blocked road.
 Xiao et al. [22] (2009)
introduced the definition of the risk of  paths on a network, which
is really a function on the set of all the paths with a same source
on the network; and introduced also the anti-risk path (ARP) problem
of finding a path such that the solution of its has minimum risk,
which, on the one hand, is a kind of single-source shortest path
problem, and on the other hand, is different from the classical
single-source shortest path problem. They  showed also the  ARP
problem can be solved in $O(mn+n^{2}log n)$ time supposed that at
most one edge may be blocked. Afterwards, Mahadeokar and Saxena [13]
(2014) proposed a faster algorithm to solve the  ARP problem.

Srivastava and Tyagi [18] (2013),  by modifying the Prim's
Algorithm, proposed an algorithm to find the shortest path in a
specific network that consists of host systems on land and
satellites in air. Murota and Shioura [16] (2014), from the
viewpoint of discrete convex analysis and linear programming
formulation, showed that the shortest path problem can be seen as a
special case of L-concave function maximization;  and solving the LP
dual of the shortest path problem with the steepest ascent algorithm
for L-concave function maximization is exactly coincident to solving
the shortest path problem with Dijkstra's algorithm. Feng [7](2014)
presented a new exact algorithm to solve $k$ shortest simple paths
(KSP) in a network; and demonstrated that the algorithm performs
significantly better than the existing exact polynomial time
algorithm that have polynomial worst-case complexity. Basing on the
concept of the replacement path and the real time detour path, Zhang
et al. [23] (2015) proposed the definition of shortest path set and
the definition of optimal shortest path set, on an undirected graph
with source node $s$ and destination node $t$; they also, under some
conditions, investigated and presented the polynomial time
algorithms to compute the optimal shortest path set. By combining
the flow shop scheduling problem and the shortest path problem, Nip
et al. [17] (2015) first developed a synthetical optimization
problem; then they discussed the complexity to compute the solutions
and separately proposed two approximation algorithms, for the case
that the number of machines is an input and for the case that the
number of machines is fixed. Meng et al. [14] (2016) reported the
results of their experiment and research on the multiple shortest
path algorithms. By combining Bellman-Ford algorithm and Dijkstra
algorithm, Dinitz and Itzhak [5] (2017) presented a new hybrid
algorithm for the single - source shortest path problem with general
edge costs, which can improve the running time bound of Bellman-Ford
algorithm for graphs with a sparse distribution of negative cost
edges; and made some related researches; in addition, they also
suggested a new straightforward proof that the Bellman-Ford
algorithm produces a shortest paths tree. Noting the wide and
important applications of Dijkstra algorithm and Bellman-Ford
algorithm in the field of computer and software sciences, and many
researches interested in the problem to solve the shortest path
problem in these applications, AbuSalim et al. [1] (2020) made a
more profound comparison between the two popular algorithms on the
complexity and the performance in terms of shortest path
optimization. Xia D W et al. [21] (2024), based Ant Colony
Optimization ACO algorithm, proposed BiA*-ACO algorithm to recommend
the fastest route for taxicabs in a complex urban road network; and
showed it is at least $49.81\%$ more efficient than the algorithms
ACO, DA and MBFA.

In addition, the types and applications of functional are expanding,
and the types and applications of networks are expanding, see the
literatures [4,12,19,20] for example.

Motivated by the stated background of researches above, in
particular the work of the literatures [22,13] (Xiao et al. 2009,
Mahadeokar and Saxena 2014), the present work, through generalizing
the total weight path function of networks as the path functional of
graphs, develops a general single-source shortest path problem
(GSSSP), which include the classical problem CSSSP and the ARP
problem  as its special cases; and tries to design an Extended
Dijkstra's Algorithm (EDA) and an Extended Moore-Bellman-Ford
Algorithm (EMBFA) to solve the problem GSSSP under certain
conditions, which respectively reduce to Dijkstra's Algorithm and
Moore-Bellman-Ford Algorithm when the problem GSSSP is the classical
problem CSSSP. Moreover, we make some related studies to analyse and
believe the two designed algorithms.

The rest of this article is organized as follows. Some preliminaries
are presented in Section 2. Section 3 formulates the problem GSSSP.
Section 4 is specially devoted to designing algorithm EDA and
algorithm EMBFA respectively. Section 5 makes the analyses of EDA
and EMBFA.
 Section 6 shows the applications of EDA and EMBFA with a few instances of problem GSSSP.
 Finally, the paper is concluded with
the Section 7.

\section{Preliminaries}
\label{sec:1}
This section provides some preliminaries for our
sequel research.

\subsection{Conceptual framework}

(1) Suppose $V$ is a set of $n (>1)$ points.

Let $u,v\in V$ and $u\neq v$. We use $[u,v]$ to denote an edge
connecting two points $u$ and $v$. And use $(u,v)$ ($(v,u)$) to
denote a road from $u$ to $v$ ($v$ to $u$) on the edge $[u,v]$ .
(Note: $[u,v]=[v,u]$, while $(u,v)\neq(v,u)$.)

When there are more than one edge between $u$ and $v$, namely there
are the parallel edges between $u$ and $v$, $[u,i,v]$ may be used to
denote the $i$th edge; and $(u,i,v)$ may be used to denote the $i$th
road from $u$ to $v$. However, to simplify in notation afterwards,
$[u,i,v]$ is denoted as $[u,v]$ and $(u,i,v)$ is denoted as $(u,v)$
when $i$ needn't be indicated.

A edge $[u,v]$ is called undirected (/directed) if there are two
roads $(u,v)$ and $(v,u)$ (/there is only one road $(u,v)$ or
$(v,u)$) on it.

For set $A$, we use $|A|$ denotes the number of all elements in the
set $A$.

Let $E$ be all the edges and $R$ be all the roads. The triple
$(V,E,R)$ is called  as a graph, which is also denoted by the tuple
$(V,E)$ (/$(V,R)$) when all the roads $R$ (/edges $E$) needn't be
indicated for clearness and briefness. Graphs without parallel edges
are called simple. For graph $(V,E,R)$, an element of $V$ is called
a vertex of the graph. In this work we always assume that $|E|$ is
finite.

A graph is called as an undirected graph (/directed graph)
 if it has only undirected edges (/has
only directed edges)

~\\
(2) Suppose $G=(V,E,R)$ is a graph.

When $[v_{(i-1)},v_{i}]\in E, i=1,2,\cdots,k$, the orderly
combination of edges
$$\{[v_{0},v_{1}],[v_{1},v_{2}],\cdots,[v_{k-1},v_{k}]\}
([v_{(i-1)},v_{i}]\neq [v_{(j-1)},v_{j}], i\neq j)$$ is called as a
chain connecting  $v_{0}$ and $v_{k}$, denoted by
$C[v_{0},v_{1},v_{2},\cdots,v_{k-1},v_{k}]$. (Note:
$C[v_{0},v_{1},v_{2},\cdots,v_{k-1},v_{k}]=C[v_{k},v_{k-1},\cdots,v_{2},v_{1},v_{0}]$.)
If $v_{0}=v_{k}$, the chain is called as a cycle (chain).

When $(v_{(i-1)},v_{i})\in R, i=1,2,\cdots,k$, the orderly
combination of roads
$$\{(v_{0},v_{1}),(v_{1},v_{2}),\cdots,(v_{k-1},v_{k})\}$$ is called
as a path from $v_{0}$ to $v_{k}$, denoted by
$P(v_{0},v_{1},v_{2},\cdots,v_{k-1},v_{k})$, in briefness denoted by
$P(v_{0},v_{k})$ if there will be no confusion; $v_{k}$ is said to
be reachable from $v_{0}$. (Note:
$P(v_{0},v_{1},v_{2},\cdots,v_{k-1},v_{k})\neq
P(v_{k},v_{k-1},\cdots,v_{2},v_{1},v_{0})$.) We also use
$(v_{0},v_{1},\cdots,v_{k-1},v_{k})$, or
$(v_{0},v_{1})+(v_{1},v_{2})+\cdots+(v_{k-1},v_{k})$, or
$(v_{0},v_{1},\cdots,v_{k-1})+(v_{k-1},v_{k})$ to denote
$P(v_{0},v_{k})$. That is,
$$\begin{array}{lcl}& &P(v_{0},v_{k})=\{(v_{0},v_{1}),(v_{1},v_{2}),\cdots,(v_{k-1},v_{k})\}=
(v_{0},v_{1},\cdots,v_{k-1},v_{k})\\
&=&(v_{0},v_{1})+(v_{1},v_{2})+\cdots+(v_{k-1},v_{k})=(v_{0},v_{1},\cdots,v_{k-1})+(v_{k-1},v_{k}).\end{array}$$
$v_{0}$ and $v_{k}$  are respectively called the source and the
terminal of path $P$, denoted by $s(P)$ and $t(P)$;
$(v_{(i-1)},v_{i}), 0<i\leq k,$ is called a road of $P$, denoted by
$(v_{(i-1)},v_{i})\in P$. For path $P=(v_{0},v_{1},\cdots,v_{k-1})$
and path $P'=(v_{0},v_{1},\cdots,v_{k-1},v_{k})=P+(v_{k-1},v_{k})$,
we call $P$ as the farther of $P'$, denoted by $F(P')$; and $P'$ as
a son of $P$, denoted by $S(P)$. A path
$P=(v_{0},v_{1},\cdots,v_{k-1},v_{k})$ is called no cycle if
$v_{i}\neq v_{j}$ while $i\neq j$.

A mapping $w:E(/R)\rightarrow (-\infty,\infty)$ is called as  a
weight of edges (/roads). And the triple $(V,w,E)$ (/$(V,w,R)$) is
called as  a network with edge (/road) weight. The tuple $(G,w)$ is
used to represent both the networks $(V,w,E)$ and $(V,w,R)$ when
$w([u,v])=w((u,v)), \forall u,v\in V$.

Let $s\in V$ (called source point). A path of graph $G$ with the
 source point $s$ is called a path of $[G,s]$. Some paths
of $[G,s]$ is called a path system on $[G,s]$. All the paths of
$[G,s]$ is called the complete path system on $[G,s]$. All the no
cycle paths of $[G,s]$ is called the no cycle path system on
$[G,s]$. To be convenient and clear, we stipulate $(s,s)$ is a
special road and a special path with the source point $s$, and
$(s,s)+(s,v)=(s,v)$ for each road $(s,v)$.

~\\
(3) Suppose $\mathcal{P}$ is a path system on $[G,s]$.

Put $\mathcal{P}(u)=\{P\in\mathcal{P}|t(P)=u\}$ (particulary,
$\mathcal{P}(s)=\{(s,s)\})$,
$V(\mathcal{P})=\{u|\mathcal{P}(u)\neq\emptyset\}$, $
R(\mathcal{P})= \{(u,v)\in P|P\in\mathcal{P}\}$ and $E(\mathcal{P})=
\{[u,v]|(u,v)\in R(\mathcal{P})\}$. Moreover, $\mathcal{P}_{nc}$
denotes all the no cycle pathes in $\mathcal{P}$, and
$\mathcal{P}_{nc}(u)$ denotes all the no cycle paths in
$\mathcal{P}(u)$.

$\forall P\in[\mathcal{P}\setminus\{(s,s)\}]$, define $(s,s)\prec
P$; $\forall P,P'\in[\mathcal{P}\setminus\{(s,s)\}]$, define $P\prec
P'$ if and only if $P=(v_{0},v_{1},\cdots,v_{k}),
P'=(v_{0},v_{1},\cdots,v_{k},v_{k+1},\cdots,v_{k'}), 0\leq k<k'$;
and $\forall P, P'\in\mathcal{P}$, define $P\preceq P'$ if and only
if $P\prec P'$ or $P=P'$.

A mapping $f:\mathcal{P}\rightarrow (-\infty, \infty)$ is called as
a path functional on $\mathcal{P}$.

Finally, $\forall v\in V(\mathcal{P})$,
 define $m_f(v)=\inf\{f(P)|P\in \mathcal{P}(v)\}$. And $P\in \mathcal{P}$ is called a
 shortest (minimum) path on $f$ if $f(P)=m_f(t(P))$.

Below, we always assume that $G$ is a graph, $s\in V(G)$ is a
source, $\mathcal{P}$ is a path system on $[G,s]$, $f$ is a path
functional on the system $\mathcal{P}$ and $f((s,s))=0$.

 \subsection{Basic definitions}

\textbf{Definition 1} (\romannumeral1) $f$ is said to be
non-decreasing if and only if $\forall P,P'\in \mathcal{P}$,
provided
 $P\preceq P'$, we have $f(P)\leq f(P')$.

(\romannumeral2) $f$ is said to be increasing if and only if
$\forall P,P'\in \mathcal{P}$, provided
 $P\prec P'$, we have $f(P)< f(P')$.

(\romannumeral3) $f$ is said to be  weak
 order-preserving
(WOP) if and only if $\forall u,v\in V$,  provided $\forall P,P'\in
\mathcal{P}(u)$,
 $P+(u,v),P'+(u,v)\in \mathcal{P}(v)$ and $f(P)<f(P')$, we
have $f(P+(u,v))<f(P'+(u,v))$.

(\romannumeral4) $f$ is said to be semi
 -order-preserving
  (SOP) if and only if
$\forall u,v\in V$, provided $\forall P,P'\in \mathcal{P}(u)$,
$P+(u,v),P'+(u,v)\in \mathcal{P}(v)$ and $f(P)\leq f(P')$, we have
$f(P+(u,v))\leq f(P'+(u,v))$.

(\romannumeral5) $f$ is said to be
 order-preserving
 (OP) if and only if $f$ is WOP, and
$\forall u,v\in V$, provided $\forall P,P'\in \mathcal{P}(u)$,
$P+(u,v),P'+(u,v)\in \mathcal{P}(v)$ and $f(P)=f(P')$, we have
$f(P+(u,v))=f(P'+(u,v))$.

~\\
\textbf{Definition 2} $f$ is said to have no negative
(/non-positive) cycle if and only if $\forall v\in V$, provided
$\forall P\in \mathcal{P}(v)$ and $\forall P'=P+
(v,v_1,\cdots,v_k,v)\in \mathcal{P}(v)$, we have $f(P')-f(P)\geq 0$
(/$f(P')-f(P)>0$). $f$ is called conservative if it has no negative
cycle.

~\\
\textbf{Definition 3} $f$ is said to be
 weak inherited on shortest path (WISP) if $\forall v\in
[V(\mathcal{P})\setminus\{s\}]$, provided that the shortest path
from $s$ to $v$ exists, then there must be a path
$P=(v_0,v_1,\cdots,v_k)\in\mathcal{P}(v), k\geq 1,$ such that
$P_i=(v_0,v_1,\cdots,v_i)$ is the shortest path, $i=1,2,\cdots,k$.
(Note: $v_k=v$).
 $f$ is said to be inherited on shortest path (ISP) if $\forall P\in [\mathcal{P}\setminus\{(s,s)\}]$, provided that $P$ is
the shortest path, then $F(P)$ must be the shortest path.
\subsection{Basic propositions}
For the above definitions, we have the following statements.

~\\
\textbf{Proposition 1} $\forall v\in V$, $|\mathcal{P}_{nc}(v)|$ is
finite.

~\\
\textbf{Proposition 2} If $f$ is increasing, then it must have no
non-positive cycle. If $f$  has no non-positive cycle  or is
non-decreasing, then it must have no negative cycle; that is, $f$ is
conservative.

~\\
\textbf{Proposition 3} If $f$ is OP, then it must be WOP and SOP.

~\\
\hspace*{0.5cm}For the proofs of the three propositions above is
trivial, here we omit them.

~\\
\textbf{Proposition 4} Let $\mathcal{P}$ be the complete path
system. (\romannumeral1) If $f$ has no non-positive cycle and is
WOP, then, $\forall v\in [V(\mathcal{P})\setminus\{s\}]$, the
shortest path from $s$ to $v$ has no cycle, and
$m_f(v)=\min\{f(P)|P\in \mathcal{P}_{nc}(v)\}$. (\romannumeral2) If
$f$ is conservative and SOP, then, $\forall v\in
[V(\mathcal{P})\setminus\{s\}]$,
 $m_f(v)=\min\{f(P)|P\in \mathcal{P}_{nc}(v)\}$. (\romannumeral3) If $f$ is
conservative, SOP and WISP,  then, $\forall v\in
[V(\mathcal{P})\setminus\{s\}]$, there is a path
$P=(v_0,v_1,\cdots,v_{k-1},v_{k})\in\mathcal{P}_{nc}(v) (k\geq 1)$
such that $P_i=(v_0,v_1,\cdots,v_i)$ is the shortest path for any
$i=0,1,2, \cdots,k$.

\textbf{Proof} $\forall v\in [V(\mathcal{P})\setminus\{s\}]$, let
$P$ be a path from $s$ to $v$ and have cycles. We can assume
$P=(v_0,v_1,\cdots,v_l,v'_1,v'_2,\cdots,v'_{l'},v_l,v_{l+1},v_{l+2},\cdots,
v_k)$ ($v_k=v$), $l'\geq 1$. Put $P_1=(v_0,v_1,\cdots,v_l)$,
$P_2=(v_0,v_1,\cdots,v_l,v'_1,v'_2,\cdots,v'_{l'},v_l)$ and
$P_3=P_1+(v_{l},v_{l+1})+\cdots+(v_{k-1},v_k)=(v_0,v_1,\cdots,v_l,v_{l+1},v_{l+2},\cdots,v_{k-1},v_k)$.
 For $\mathcal{P}$ is the complete path system,
 $P_1,P_2,P_3\in\mathcal{P}$.

(1) For $f$ has no non-positive cycle, we have $f(P_2)>f(P_1)$.
Since also $f$ is WOP, we further have
$f(P_2+(v_{l},v_{l+1}))>f(P_1+(v_{l},v_{l+1})),\cdots,f(P_2+(v_{l},v_{l+1})
+\cdots+(v_{k-1},v_k))>f(P_1+(v_{l},v_{l+1})+\cdots+(v_{k-1},v_k))$.
Thus $f(P)>f(P_3)$. This implies that $P$ can not be a shortest path
from $s$ to $v$ and leads to $m_f(v)=\inf\{f(P)|P\in
\mathcal{P}(v)\}\leq\min\{f(P)|P\in
\mathcal{P}_{nc}(v)\}\leq\inf\{f(P)|P\in
\mathcal{P}(v)\}=m_f(v)\Rightarrow m_f(v)=\min\{f(P)|P\in
\mathcal{P}_{nc}(v)\}$. Hence (\romannumeral1) holds.

(2) For $f$ is conservative, we have $f(P_2)\geq f(P_1)$. Since also
$f$ is SOP, we further have $f(P_2+(v_{l},v_{l+1}))\geq
f(P_1+(v_{l},v_{l+1})),\cdots,f(P_2+(v_{l},v_{l+1})
+\cdots+(v_{k-1},v_k))\geq
f(P_1+(v_{l},v_{l+1})+\cdots+(v_{k-1},v_k))$. Thus $f(P)\geq
f(P_3)$. This implies that there must be a shortest path from $s$ to
$v$ such that it has no cycle. Hence (\romannumeral2) holds.

(3) Let $v\in [V(\mathcal{P})\setminus\{s\}]$. In terms of the
conclusion (\romannumeral2), there must be a shortest path from $s$
to $v$. Since also $f$ is WISP, we can further know there is a path
$(v_0,v_1,v_2,\cdots,v_k)\in \mathcal{P}(v)$ such that
$(v_0,v_1,v_2,\cdots,v_i), i=0,1,2,\cdots,k$, are all the shortest
pathes. Finally, following the approach to prove conclusion
(\romannumeral2), we can easily prove that there is  a  path
$P=(v_0,v'_1,v'_2,\cdots,v'_{k'})\in \mathcal{P}_{nc}(v)$ such that
$(v_0,v'_1,v'_2,\cdots,v'_i), i=0,1,2,\cdots,k'$, are all the
shortest paths. Hence (\romannumeral3) holds.

The proof completes. $\hfill\square$

~\\
\textbf{Corollary 1} Let $\mathcal{P}$ be the complete path system.
If $f$ is non-decreasing and SOP, then $\forall v\in V(\mathcal{P}),
m_f(v)=\min\{f(P)|P\in \mathcal{P}_{nc}(v)\}$.

~\\
\textbf{Lemma 1} Let $\mathcal{P}$ be the complete path system.
Suppose also $f$ has no non-positive cycle. If the path
$P=(v_0,v_1,\cdots,v_k,v_{(k+1)},v_{(k+2)},\cdots,v_{(k+l)})$
($k\geq 0, l\geq 1$) such that
$$P_{(k+i)}=(v_0,v_1,\cdots,
v_k,v_{(k+1)},v_{(k+2)},\cdots, v_{(k+i)}), i=1,2,\cdots,l,$$ are
all the shortest paths, then $v_{(k+i)}\neq v_{(k+j)}$ while $1\leq
i<j\leq l$.

 \textbf{Proof} It is obvious that $v_{(k+i)}\neq v_{(k+j)}$ when
$j=i+1$. Assume that $1\leq i<j\leq l, j-i\geq 2,$ and $v_{(k+i)}=
v_{(k+j)}=v$. Then, $f(P_{(k+i)})=f(P_{(k+j)})=m_f(v)$. For $f$ has
no non-positive cycle, we have $f(P_{(k+j)})>f(P_{(k+i)})$. This is
contradictory to $f(P_{(k+i)})=f(P_{(k+j)})$. Hence the lemma holds.
$\hfill\square$

~\\
\textbf{Proposition 5} Let $\mathcal{P}$ be the complete path
system. If $f$ has no non-positive cycle and is SOP, then $f$ is
WISP.

\textbf{Proof} $\forall v\in [V(\mathcal{P})\setminus\{s\}]$, let
$P=(v_0,v_1,\cdots,v_k,v), k\geq 0,$ be a shortest path. If
$(v_0,v_1,\cdots,v_k)$ is not the shortest path, then, from the term
(\romannumeral2) of Proposition 4, there must be a shortest path
$P'=(v_0,v'_1,\cdots,v'_{k'},v_{k})$. For $\mathcal{P}$ is complete,
$[P'+(v_{k},v)]\in \mathcal{P}$. Since $f$ is SOP, we have
$f(P'+(v_{k},v))\leq f(P)$. This implies that $[P'+(v_{k},v)]$ is
also the shortest path. If $(v_0,v'_1,\cdots,v'_{k'})$ is not the
shortest path, then, in the same way, there must be a shortest path
$P''=(v_0,v''_1,\cdots,v''_{k''}, v'_{k'})$ such that
$[P''+(v'_{k'},v_{k})], [P''+(v'_{k'},v_{k})+(v_{k},v)]$ are all the
shortest path. $\cdots$ Assume that the step has been performed $l$
times. Then we can obtain a path $P^\ast=(v_0,v^\ast_1,\cdots,
v^\ast_{k^\ast},v^\ast_{(k^\ast+1)},v^\ast_{(k^\ast+2)},\cdots,v^\ast_{(k^\ast+l)})\in\mathcal{P}(v)$
($k^\ast\geq 0, l\geq 1$) such that
$$P_{(k^\ast+i)}=(v_0,v^\ast_1,\cdots,
v^\ast_{k^\ast},v^\ast_{(k^\ast+1)},v^\ast_{(k^\ast+2)},\cdots,
v^\ast_{(k^\ast+i)}), i=1,2,\cdots,l,$$ are all the shortest paths.
By Lemma 1,  we have $v^\ast_{(k^\ast+i)}\neq v^\ast_{(k^\ast+j)}$
while $1\leq i<j\leq l$. This implies that $l\leq n$. Therefore, we
can eventually get a shortest path
$\bar{P}=(v_0,\bar{v}_1,\cdots,\bar{v}_{\bar{k}},v)$ such that
$$\bar{P}_i=(v_0,\bar{v}_1,\cdots,
\bar{v}_{\bar{k}}), i=1,2,\cdots,\bar{k},$$ are all the shortest
paths. Hence the proposition holds. $\hfill\square$

~\\
\textbf{Proposition 6} Let $\mathcal{P}$ be the complete path
system. If $f$ is conservative, SOP and WOP, then $f$ is ISP.

\textbf{Proof} Let $P=(v_0,v_1,\cdots,v_k)\in\mathcal{P},k\geq 1$,
be a shortest path. If $F(P)$ is not the shortest path, then, from
the term (\romannumeral2) of Proposition 4, there must be another
path $P'\in\mathcal{P}_{nc}(v_{k-1})$ such that
$f(P')=m_f(v_{k-1})<f(F(P))$. For $\mathcal{P}$ is complete path
system, we have $[P'+(v_{k-1},v_{k})]\in \mathcal{P}(v_{k})$. Since
also $f$ is WOP and $[P'+(v_{k-1},v_{k})],[F(P)+(v_{k-1},v_{k})]\in
\mathcal{P}(v_{k})$, we further have
$f(P'+(v_{k-1},v_{k}))<f(F(P)+(v_{k-1},v_{k}))=f(P)$, which is
contradictory with that $P$ is  the shortest path. Hence the
proposition holds. $\hfill\square$

~\\
\textbf{Corollary 2} Let $\mathcal{P}$ be the complete path system.
If $f$ is conservative and OP, then $\forall v\in V(\mathcal{P}),
m_f(v)=\min\{f(P)|P\in \mathcal{P}_{nc}(v)\}$  and $f$ is ISP.

\textbf{Proof} From Proposition 4 and Proposition 6, we can easily
know the corollary holds. $\hfill\square$

~\\
\textbf{Remark 1} The propositions above are not only the basis for
our to design and study the next algorithms EDA and EMBFA, but also
fully shows that the contents of the path functional, especially its
order relations, are very profound, extensive and interesting, which
means that for the path functional,
 there are still many problems needing to
research.

\section{Problem}
\label{sec:2}

\textbf{Definition 4} The problem to find a path
$P\in\mathcal{P}(v)$ such that $f(P)=m_{f}(v)$ for all $v\in
V(\mathcal{P})$ is called as general single-source shortest path
problem (GSSSP) on $[G,s,\mathcal{P},f]$.

~\\
\hspace*{0.5cm}It is clear that the problem GSSSP is just the
problem CSSSP when $G$ is the graph with weight $w$ and $f$ is the
path function $d$ in the example 1 of Section 6. It is also clear
that the  ARP problem, see Xiao et al. [22] (2009) or example 2, is
an instance of the problem GSSSP. The two facts show that  the
problem GSSSP is really generalization of the problem CSSSP.

~\\
\textbf{Theorem 1} Let $\mathcal{P}$ be the complete path system. If
$f$ is conservative and SOP, then the problem GSSSP can be solved.
That is, $\forall v\in V(\mathcal{P})$, there is a path
$P\in\mathcal{P}(v)$ such that $f(P)=m_{f}(v)$, namely $P$ is a
shortest path from $s$ to $v$.

\textbf{Proof}  From the term (\romannumeral2) of Proposition 4, the
theorem 1 holds. $\hfill\square$

\section{Algorithms}
\label{sec:3} For  problem GSSSP, following the approaches of the
algorithms DA and MBFA, an extended Dijkstra's algorithm (EDA) and
an extended Moore-Bellman-Ford algorithm (EMBFA) can be respectively
designed to solve it under certain conditions. We accomplish the
tasks in this section.

~\\
\textbf{Extended Dijkstra's Algorithm} (EDA)

\textbf{Input:} graph $G=(V,E,R)$ and vertex $s\in V$, with a path
functional $f$ on the complete path system $\mathcal{P}$ of $[G,s]$,
which is non-decreasing and SOP.

\textbf{Output:} a path system $\mathcal{T}$ on $[G,s]$ and the
graph $T=(V(\mathcal{T}),E(\mathcal{T}),R(\mathcal{T}))$.

\textbf{Process:}

1. Put $P_{T}[s]=(s,s),f(P_{T}[s])=f((s,s))=0$;

\quad $\forall v\in [V\setminus\{s\}]$, set
$P_{T}[v]\leftarrow(s,\infty,v), f(P_{T}[v])\leftarrow +\infty$.

\quad Set $C\leftarrow\emptyset$, $\mathcal{T}\leftarrow\emptyset$.
(Set $k\leftarrow (-1)$.)

2.

(1) Find a $u\in[V\setminus C]$ such
that$f(P_{T}[u])=\min\{f(P_{T}[v])|v\in [V\setminus C]\}$.

(2) Set $C\leftarrow [C\bigcup\{u\}]$,
$\mathcal{T}\leftarrow[\mathcal{T}\bigcup\{P_{T}[u]\}]$.

\quad\quad$\forall v\in [V\setminus C]$, if $ (u,v)\in \mathcal{P}$
and $f(P_{T}[v])>f(P_{T}[u]+(u,v))$, set

\quad\quad$P_{T}[v]\leftarrow(P_{T}[u]+(u,v))$.

\quad\quad(Set $k\leftarrow (k+1)$, then put $v(k)=u$.)

(3) If $[V\setminus C]=\emptyset$ or $\min\{f(P_{T}[v])|v\in
[V\setminus C]\}=+\infty$, go to the step 3.

\quad\quad Otherwise, return to step 2.

3. Output the path system $\mathcal{T}$ and the graph
$T=(V(\mathcal{T}),E(\mathcal{T}),R(\mathcal{T}))$.

\quad (Output the vertices: $v(i), i=0,1,2,\cdots,
|V(\mathcal{T})|-1$.) Then stop.

~\\
\textbf{Remark 2} To be concise, please first understand the EDA and
the related proof (see next section) under the condition that $G$ is
a simple graph. When $G$ is not a simple graph, we should consider
the parallel edges. For instance, when there exist two parallel
roads $(u,1,v), (u,2,v)\in R(\mathcal{P})$ and $u$ is $v(k)$, on the
update of $P_{T}[v]$ in the $(k+1)$ time iteration,
$P_{T}[u]+(u,1,v)$ and $P_{T}[u]+(u,2,v)$ should be  simultaneously
involved in the comparison and the replacement. More specifically,
we should understand and execute the term (2) of step 2 as follows.

Set $C\leftarrow [C\bigcup\{u\}]$,
$\mathcal{T}\leftarrow[\mathcal{T}\bigcup\{P_{T}[u]\}]$.

$\forall v\in [V\setminus C]$, if $ (u,i,v)\in \mathcal{P},
i=1,2,\cdots, |I(u,v)|$, $I(u,v)=\{$all the roads $(u,i,v)\}$, and
$f(P_{T}[v])>\min\{f(P_{T}[u]+(u,i,v))|i=1,2,\cdots, |I(u,v)|\}$,
find $i'$ such that
$f(P_{T}[u]+(u,i',v))=\min\{f(P_{T}[u]+(u,i,v))|i=1,2,\cdots,
|I(u,v)|\}$, set

\quad\quad $P_{T}[v]\leftarrow(P_{T}[u]+(u,i',v))$.

\quad\quad(Set $k\leftarrow (k+1)$, then put $v(k)=u$.)

~\\\hspace*{0.5cm} In order to simplify the analytical process of
algorithm EDA, we also propose the next algorithm STA.

~\\
\textbf{Spanning Tree Algorithm }(STA)

\textbf{Input:} graph $G=(V,E,R)$, point $s\in V$.

\textbf{Output:} a graph $T$.

\textbf{Process:} 1. Set  $C\leftarrow\{s\}, R(T)\leftarrow
\{(s,s)\}$.

\quad\quad \quad\qquad(Set $k\leftarrow0$, then put $v(k)=s$,
$T_k=(C,R(T))$.)

If $\{(u,v)|u\in C, v\in[V\setminus C], (u,v)\in R\}\neq\emptyset$,
implement the next step. Otherwise, go to step 3.

2. Find a $u\in C$ and a $v\in[V\setminus C]$  such that $(u,v)\in
R$. Then set
 $$C\leftarrow [C\cup\{v\}],R(T)\leftarrow [R(T)\cup\{(u,v)\}].$$
\qquad \qquad\qquad\qquad(Set $k\leftarrow (k+1), v(k)=v$,
$T_k=(C,R(T))$.)

If $\{(u,v)|u\in C, v\in[V\setminus C], (u,v)\in
R(\mathcal{P})\}\neq\emptyset$, return to step 2. Otherwise,
implement the next step.

3.  Put $V(T)=C$ and $T=(V(T),R(T))$. Output the graph $T$. Then
stop.

~\\
\textbf{Extended Moore-Bellman-Ford Algorithm} (EMBFA)

\textbf{Input:} graph $G=(V,E,R)$ and vertex $s\in V$, with a path
functional $f$ on the complete path system $\mathcal{P}$ of $[G,s]$,
which is conservative and OP.

\textbf{Output:} a path system $\mathcal{T}$ on $[G,s]$ and the
graph $T=(V(\mathcal{T}),E(\mathcal{T}),R(\mathcal{T}))$.

\textbf{Process:}

1.

Put $P_{T}[s]=(s,s),f(P_{T}[s])=0$;

\quad set $P_{T}[v]\leftarrow(s,\infty,v), f(P_{T}[v])\leftarrow
+\infty, \forall v\in [V\setminus\{s\}]$.

\quad Set $\mathcal{T}\leftarrow\{P_{T}[s]\}$.

2. For $i=1,2,\cdots,n$, do:

\quad for each road $(u,v)\in [R(\mathcal{P})\setminus\{(s,s)\}]$,
if $P_{T}[u]\in \mathcal{T}$ and $$f(P_{T}[v])>f(P_{T}[u]+(u,v)),$$

\quad then in turn set: $\mathcal{T'}\leftarrow[\mathcal{T}\setminus
\{P_{T}[v]\}]$; $P_{T}[v]\leftarrow P_{T}[u]+(u,v)$;

\quad and $\mathcal{T}\leftarrow[\mathcal{T'}\cup \{P_{T}[v]\}]$.

\quad  (Set $T=(V(\mathcal{T}),E(\mathcal{T}),R(\mathcal{T}))$.)

3. Output the path system $\mathcal{T}$ and the graph
$T=(V(\mathcal{T}),E(\mathcal{T}),R(\mathcal{T}))$.

\quad Then stop.

~\\
\textbf{Remark 3} For EDA, $\forall v\in V(\mathcal{P})$, $P_{T}[v]$
will not change after it become a member of $\mathcal{T}$. However,
for EMBFA, some $P_{T}[v]$ may change after it become a member of
$\mathcal{T}$. The two facts are useful for us to understand the
algorithms  EDA and EMBFA.

\section{Analysis of algorithms}
\label{sec:4} \textbf{Lemma 2} For algorithm STA, we have the
following conclusions.

(\romannumeral1)  The algorithm  works well. (That is, it can
effectively input, output and stop.) (\romannumeral2) The output $T$
of STA  is a spanning tree of $\bar{G}$. (\romannumeral3) $T$ is an
arborescence rooted at $s$. (\romannumeral4) The running time of STA
is no more than $O(mn)$. Here
$\bar{G}=(V(\mathcal{P}),E(\mathcal{P}),R(\mathcal{P}))$,
$\mathcal{P}$ is the complete path system on $[G,s]$ and $m=|E(V)|$.

\textbf{Proof} Obviously, we always have that $1\leq |C|\leq n=|V|$
and $|C|$ is increasing. This implies that the step 2 is performed
no more than n times. Hence (\romannumeral1) holds.

In terms of the process of Algorithm STA and (\romannumeral1), it is
obvious that $T$ is connected subgraph of $\bar{G}$. Note that
$\{(u,v)|u\in C, v\in[V\setminus C], (u,v)\in R\}\neq\emptyset$ so
long as $|C|<|V(\bar{G})|$. We can also know that
$|V(T)|=|V(\bar{G})|$. Therefore, in order to show (\romannumeral2)
is true, we
 need only to prove that $T$ has no cycle. Clearly, $T_{0}$ has no
 cycle. Suppose that $T_{k}$ has no
 cycle. Then, for $v(k+1)\neq v(0), v(1), \cdots,  v(k)$, $T_{(k+1)}$ also has no
 cycle. This implies that $T$ has no cycle. Hence
(\romannumeral2) holds.

By the process of STA, we can easily know that for any edge $[u,v]$
of $T$, there exists only one road ($(u,v)$ or $(v,u)$), that is,
$T$ is a directed graph. Moreover, we can easily know also that
$|\delta^-(s)|=0$ and $|\delta^-(v)|=1, \forall v\in
[V(T)\setminus\{s\}]$. Here $|\delta^-(v)|$ is the in-degree of $v$
in the directed graph $T$.  Therefore $T$ is an arborescence rooted
at $s$. Hence (\romannumeral3) holds.

Finally, it is obvious that the running time of STA depends on the
complexity of the step 2. However, the number of total iterations
 of the
 step 2 do not exceed $n$ and the complexity of each run of the
 step 2 does
not exceed $O(m)$. Hence (\romannumeral4) holds.

The proof completes. $\hfill\square$

~\\
\textbf{Remark 4} For the necessary basic concepts of graphs, such
as  spanning tree,
  arborescence, induced subgraph,
the degree of a vertex, the in-degree of a vertex, so on, please see
2.1 and 2.2 of the monograph [11] (Korte and Vygen 2000) or other
related books.

~\\
\textbf{Remark 5} For algorithm STA, by properly designing the
method to find the road $(v,u)$ in step 2, the running time can be
greatly simplified.  Please, see the Graph Scanning Algorithm in 2.3
of the monograph (Korte and Vygen 2000) for reference. However, the
current STA can more effectively implement its functionality in the
present work, which is helpfully to prove the following theorem 2,
so we propose it in the current pattern.

~\\
\textbf{Theorem 2} For algorithm EDA, we have the following
conclusions.

(\romannumeral1) The algorithm works well. (\romannumeral2) The
output $T$ is a spanning tree of graph $\bar{G}$. (\romannumeral3)
$T$ is an arborescence rooted at $s$. (\romannumeral4) $\forall v\in
V(\mathcal{P})$, the path from $s$ to $v$ on the tree $T$  is
$P_{T}[v]$. (\romannumeral5) $\forall v\in V(\mathcal{P}),
 f(P_{T}[v])=m_f(v)$. (\romannumeral6) The running time is
$M(n)O(n^2)$, provided $\Delta(G)=O(n)$  and $f(P)$ can be obtained
in $M(n)$ time of calculation for any $P\in \mathcal{P}$.

 Here,
$\bar{G}=(V(\mathcal{P}),E(\mathcal{P}),R(\mathcal{P}))$, which is
the subgraph of $G$ induced by $V(\mathcal{P})$; $n=|V|$;
$\Delta(G)=\max\{|\delta(v)||v\in V\}$, $|\delta(v)|$ is the degree
of vertex $v$, which is called as the maximum degree of graph $G$;
$M(n)$ is a constant related to $n$.

\textbf{Proof.}  Note that $V$ and $R$ are all the finite sets. By
observing the process of algorithm EDA, we can easily know that
(\romannumeral1) holds. Note that $[V\setminus C]=\emptyset$ or
$\min\{f(P_{T}[v])|v\in [V\setminus C]\}=+\infty$ is equivalent to
$\{(u,v)|u\in C, v\in[V\setminus C], (u,v)\in
R(\mathcal{P})\}=\emptyset$. By carefully examining the second step
of EDA and the second step of STA, it can be easily known that the
second step of EDA is a specific implementation of the second step
of STA. Hence (\romannumeral2) and (\romannumeral3) can be
immediately obtained from Lemma 2. Note $P_{T}[v]\in \mathcal{T}$.
According to the relation between $T$ and $\mathcal{T}$, we can
easily know that $\mathcal{T}$ is actually the complete path system
on  $[T,s]$. Hence (\romannumeral4) holds. As $\Delta(G)=O(n)$ and
$f(P)$ can be obtained in $M(n)$ times of calculation for any $P\in
\mathcal{P}$, we can easily know that the running time of the step 2
of EDA is $M(n)O(n)$. Note that the total iterations of step 2 is no
more than $n$ time. We can further know the running time of EDA is
$M(n)O(n^2)$. That is, (\romannumeral6) holds. Next we focus to
prove (\romannumeral5).

Note that the outputted vertices:  $v(l),
l=0,1,\cdots,|V(\mathcal{P})|-1$, actually are all the elements of
$V(\mathcal{P})$. The term (\romannumeral5) can be interpreted as:
$f(P_{T}[v(l)])=m_f(v(l)), l=0,1,\cdots,|V(\mathcal{P})|-1$. We use
mathematical induction to complete the proof.

In the first place, it is obvious that (\romannumeral5) holds when
$l=0$, namely $f(P_{T}[v(0)])\\=m_f(v(0))$. In fact, we can easily
know that $v(0)=s$. Hence $P_{T}[v(0)]=(s,s)$ and
$f(P_{T}[v(0)])=f((s,s))=0=m_f(v(0))$. Assume $|\mathcal{P}|\geq2$.
Then it is also obvious that (\romannumeral5) is true when $l=1$,
namely $f(P_{T}[v(1)])=m_f(v(1))$.

Suppose (\romannumeral5) is true for  $0\leq l\leq
k<|V(\mathcal{P})|-1$ with $k\geq 1$. Then we can  prove that it is
also true for $0\leq l\leq(k+1)$.

In fact, we only need to prove that (\romannumeral5) holds for
$(k+1)$.

Assume that (\romannumeral5) does not holds for $(k+1)$, namely,
$m_f(v(k+1))<f(P_{T}[v(k+1)])$. We can show that the assumption is
not true as follows.

Let $v(k+1)=v^\ast$ and $P_{T}[v(k+1)]=(s,\cdots, \tilde{v},
v^\ast)$. Then $m_f(v^\ast)<f(P_{T}[v^\ast])$, and from term
(\romannumeral4), we have $P_{T}[\tilde{v}]=(s,\cdots, \tilde{v})$.

In terms of the assumption and Corollary 1, there is a path
$$P=(v_{0},v_{1},\cdots,v_{i},v^\ast)\in\mathcal{P}_{nc}(v^\ast)\eqno(3)$$
such that $f(P)=m_f(v^\ast)<f(P_{T}[v^\ast])$ and
$(v_{0},v_{1},\cdots,v_{i})\neq P_{T}[\tilde{v}]$.
 For the vertex
$v_{i}$, there are only three situations: (1) $
v_{i}\in\{v(l)|l=0,1,\cdots,k\}$, and $ v_{i}=\tilde{v}$; (2) $
v_{i}\in\{v(l)|l=0,1,\cdots,k\}$, but $ v_{i}\neq \tilde{v}$; (3) $
v_{i}\notin\{v(l)|l=0,1,\cdots,k\}$.

In the situation (1), $f(P_{T}[v_{i}])=m_f(v_{i})$ for
$v_{i}\in\{v(l)|l=0,1,\cdots,k\}$, and then
$f((v_{0},v_{1},\cdots,v_{i}))\geq m_f(v_{i})
=f(P_{T}[v_{i}])=f(P_{T}[\tilde{v}])$. Because $f$ is SOP, we have
$m_f(v^\ast)=f((v_{0},v_{1},\cdots,v_{i},v^\ast))=f((v_{0},v_{1},\cdots,v_{i})+(v_{i},v^\ast))\geq
f(P_{T}[\tilde{v}]+(v_{i},v^\ast))=f(P_{T}[v^\ast])$. This is in
contradiction with $m_f(v^\ast)<f(P_{T}[v^\ast])$.

In the situation (2), $f(P_{T}[v_{i}])=m_f(v_{i})\leq
f((v_{0},v_{1},\cdots,v_{i}))$. Because $f$ is SOP,
$$\begin{array}{lcl}& & f(P_{T}[v_{i}]+(v_{i},v^\ast))\\&\leq &
f((v_{0},v_{1},\cdots,v_{i})+(v_{i},v^\ast))=f((v_{0},v_{1},\cdots,v_{i},v^\ast))<
f(P_{T}[v^\ast]).\end{array}\eqno(4)$$ For $v_{i}\neq \tilde{v}$, we
have
$$P_{T}[v_{i}]+(v_{i},v^\ast)\neq
P_{T}[\tilde{v}]+(v_{i},v^\ast)=(s,\cdots, \tilde{v},
v^\ast).\eqno(5)$$ Since $v_{i}\in\{v(l)|l=0,1,\cdots,k\}$ and
$v^\ast\notin\{v(l)|l=0,1,\cdots,k\}$, on the basis of the process
of EDA, from formulae (4) and (5), we have
$P_{T}[v(k+1)]\neq(s,\cdots, \tilde{v}, v^\ast)$. This is in
contradiction with $P_{T}[v(k+1)]=(s,\cdots, \tilde{v}, v^\ast)$.

In the situation (3), for $v_{0}=s$, we can find a vertex
$$v_{i'}\in[\{v_{0},v_{1},\cdots,v_{i-1}\}\cap\{v(l)| l=0,1,\cdots,k\}]$$
such that $v_{i'+1}, v_{i'+2},\cdots,v_{i}\notin
\{v(l)|l=0,1,\cdots,k\}$. Following the approach of situation (2),
we can obtain
$$\begin{array}{lcl}& &f(P_{T}[v_{i'}]+(v_{i'},v_{i'+1}))\\&\leq &
f((v_{0},v_{1},\cdots,v_{i'})+(v_{i'},v_{i'+1}))=f(v_{0},v_{1},\cdots,v_{i'},v_{i'+1}).\end{array}\eqno(6)$$
On the other hand, for $f$ is non-decreasing, we have
$$f(v_{0},v_{1},\cdots,v_{i'},v_{i'+1})\leq f(P)<f(P_{T}[v^\ast]). \eqno(7)$$
From formulae (6) and (7), we can obtain
$f(P_{T}[v_{i'}]+(v_{i'},v_{i'+1}))<f(P_{T}[v^\ast])$. By formulae
(3), we easily know $v_{i'+1}\neq v^\ast$.  Noting
$v_{i'}\in\{v(l)|l=0,1,\cdots,k\}$, on the basis of the process of
EDA, we can derive $v(k+1)\neq v^\ast$. This is in contradiction
with $v(k+1)=v^\ast$.

Combined with the above results of the three situations, we know
that the assumption $m_f(v(k+1))<f(P_{T}[v(k+1)])$ is incorrect.
That is, $f(P_{T}[v(k+1)])=m_f(v(k+1))$. Finally, by the induction
principle, $f(P_{T}[v(l)])=m_f(v(l))$ for
$l=0,1,\cdots,|V(\mathcal{P})|-1$. Hence (\romannumeral5) holds.

The proof is completed. $\hfill\square$

~\\
\textbf{Theorem 3} For algorithm EMBFA, we have the following
conclusions.

(\romannumeral1) The algorithm EMBFA works well. (\romannumeral2)
The output
 $T$ is a spanning tree of graph $\bar{G}$. (\romannumeral3) $T$ is an arborescence rooted at
$s$. (\romannumeral4) $\forall v\in V(\mathcal{P})$,  the path from
$s$ to $v$ on the tree $T$  is $P_{T}[v]$, and $f(P_{T}[v])=m_f(v)$.
(\romannumeral5) The running time is $M(n)O(nm)$, provided $f(P)$
can be obtained in the $M(n)$ time of calculation for any $P\in
\mathcal{P}$.

Here, $\bar{G}=(V(\mathcal{P}),E(\mathcal{P}),R(\mathcal{P}))$,
which is the subgraph of $G$ induced by $V(\mathcal{P})$; $n=|V|,
m=|E|$; $M(n)$ is a constant related to $n$.

\textbf{Proof.} Note that $V$ and $R$ are all the finite sets. We
can easily know (\romannumeral1) holds. Note that $f(P)$ can be
obtained in the $M(n)$ time of calculation. Following the approach
to analyse the running time of algorithm MBFA, we can easily
 know that
(\romannumeral5) holds by the process of algorithm EMBFA. Next we
focus to prove (\romannumeral2), (\romannumeral3) and
(\romannumeral4).

For they evidently hold when $|V(\mathcal{P})|\leq 2$, we prove them
under the condition $|V(\mathcal{P})|>2$.

Above all, we claim: \\
(a) $\forall v\in V(\mathcal{P})$,  $f(P_{T}[v])=m_f(v)$; further,
once $f(P_{T}[v])$ attains $m_f(v)$, $P_{T}[v]$ \hspace*{0.5cm} will
remain unchanged in the after process; and
$P_{T}[v]$ is unique. \\
(b) $\forall v\in V(\mathcal{P})$, let
$P_{T}[v]=(v_0,v_1,\cdots,v_k, v), 0\leq k$, then
$$P_{T}[v_i]=(v_0,v_1,\cdots,v_i), i=0,1,2,\cdots,k;$$ \hspace*{0.5cm}further,
$P_{T}[v]\in \mathcal{P}_{nc}(v)$; and $k\leq |V(\mathcal{P})|-2$.
\\(c) $T$ has no cycle;
further, $T$ is directed graph; and
$$|\delta^-(s)|=0, \forall v\in
[V(\mathcal{P})\setminus\{s\}], |\delta^-(v)|=1.$$

Prove (a). When $v=s$, (a) clearly holds. For $\mathcal{P}$ is the
complete path system on $[G,s]$, $\forall v\in
[V(\mathcal{P})\setminus\{s\}]$, from Corollary 2, there is a path
$P=(v_0,v_1,\cdots,v_k,v)\in \mathcal{P}_{nc}(v), 0\leq k\leq n-2,$
such that $f(P)=m_f(v)$. From the process of EMBFA, this implies
that $P_{T}[v]$ must become a path and $f(P_{T}[v])$ must attain
$m_f(v)$ within $(k+2)$  iterations of step 2. Since also
$f(P_{T}[v])$ never increases in the process of EMBFA, we can easily
know that $f(P_{T}[v])=m_f(v)$ in the end. Thus, the first statement
is true. The correctness of the  latter two statements is obvious.
Hence (a) holds.

Prove (b). $\forall v\in V(\mathcal{P})$, let
$P_{T}[v]=(v_0,v_1,\cdots,v_k, v)$. Assume that
$P_{T}[v_i]=(v_0,v_1,\cdots,v_l,v_{l+1},\cdots,v_i), l< i\leq k$ and
$P_{T}[v_l]=(v_0,v'_1,\cdots,v'_{l'},v_l)\neq(v_0,v_1,\\\cdots,v_l)$.
Then it is implemented first that
$P_{T}[v_l]=(v_0,v'_1,\cdots,v'_{l'},v_l)$ and
$f(P_{T}[v_l])=m_f(v_l)$. By Corollary 2, $f$ is ISP, so
$f((v_0,v_1,\cdots,v_l))=m_f(v_l)=f(P_{T}[v_l])$. This implies that
$(v_0,v_1,\cdots,v_l,v_{l+1},\cdots,v_k, v)$ can no longer be the
path of $T$ (element of $\mathcal{T}$) at any stage of EMBFA, which
leads to contradictions. Thus the first statement is true.  $\forall
v\in V(\mathcal{P})$, assume that $P_{T}[v]$ has a cycle. Then we
have the following representation,
$$P_{T}[v]=(v_0,v_1,\cdots,v_l,\bar{v},v'_{1},\cdots,v'_{l'},\bar{v},\cdots,v), 1\leq l'.$$
In terms of the first statement, we have
$P_{T}[\bar{v}]=(v_0,v_1,\cdots,v_l,\bar{v})$
 and
 $P_{T}[\bar{v}]=(v_0,v_1,\cdots,v_l,\bar{v},v'_{1},\cdots,v'_{l'},\bar{v})$, which means $P_{T}[v]$ is not unique, and contradicts the claim (a).
So, $P_{T}[v]$ has no cycle, and the third statement is also true.
Hence (b) holds.

Prove (c). According to the relation between $T$ and $\mathcal{T}$,
we can easily know that if $T$ has cycles, one of the following
facts (1) and (2) holds. (1) There is a path of $\mathcal{T}$ such
that it has cycles. (2) There are two paths
$P_{T}[v],P_{T}[v']\in\mathcal{T}$ such that
$$P_{T}[v]=(v_0,v_1,\cdots,v_l,v_{(l+1)},\cdots,v_{(l+k)},\cdots,v), 2\leq k;$$
$$P_{T}[v']=(v_0,v_1,\cdots,v_l,v'_{1},\cdots,v'_{k'},\cdots,v'), 1\leq k';$$
$$v_{(l+k)}=v'_{k'}=\bar{v}, \{v'_{1},\cdots,v'_{k'-1}\}\cap\{v_{(l+1)},\cdots,v_{(l+k-1)}\}=\emptyset.$$
That is, the two paths intersect again after their separating. The
fact (1) contradicts the second statement of claim (b). By the first
statement of claim (b), the fact (2) results in
$P_{T}[\bar{v}]=(v_0,v_1,\cdots,v_l,v_{(l+1)},\cdots,v_{(l+k)})$ and
$P_{T}[\bar{v}]=(v_0,v_1,\cdots,v_l,v'_{1},\cdots,v'_{k'})$, which
means that $P_{T}[\bar{v}]$ is not unique for
$$(v_0,v_1,\cdots,v_l,v_{(l+1)},\cdots,v_{(l+k)})\neq
(v_0,v_1,\cdots,v_l,v'_{1},\cdots,v'_{k'}),$$ and contradicts the
claim (a). Thus $T$ has no cycle. Assume $T$ is not directed graph.
Then there must be $u,v\in V(T)$ such that $(u,v),(v,u)\in R(T)$.
From the approach that graph $T$ is created by path system
$\mathcal{T}$, there must be $u',v'\in V(\mathcal{T})$ such that
$(u,v)\in P_{T}[u'], (v,u)\in P_{T}[v']$, namely,
$P_{T}[u']=(s,\cdots,u,v,\cdots,u'),
P_{T}[v']=(s,\cdots,v,u,\cdots,v')$. In terms of
$P_{T}[u']=(s,\cdots,u,v,\cdots,u')$ and (b),
$P_{T}[u]=(s,\cdots,u)$ and $v\notin V(P_{T}[u])$ for $P_{T}[u']$
has no cycle. In terms of $P_{T}[v']=(s,\cdots,v,u,\cdots,v')$ and
(b), $P_{T}[u]=(s,\cdots,v,u)$ and $v\in V(P_{T}[u])$. For $v\notin
V(P_{T}[u])$ contradicts $v\in V(P_{T}[u])$, the assumption is not
true. That is,  $T$ is directed graph. For, in the beginning,
$P_{T}[s]=(s,s), f(P_{T}[s])=m_f(s)$, which remains unchanged in the
after process, it is obvious that $|\delta^-(s)|=0$. $\forall v\in
[V(\mathcal{P})\setminus\{s\}]$, we have clearly $|\delta^-(v)|\geq
1$. On the other hand, for all the paths of $T$ have the same source
point $s$ and $T$ has no cycle, we have $|\delta^-(v)|\leq 1$. So,
$|\delta^-(v)|= 1$. To sum up, (c) holds.

Finally, we show that (\romannumeral2), (\romannumeral3) and
(\romannumeral4) hold basing on the  above claims.

In terms of (a), $V(\mathcal{P})\subseteq V(\mathcal{T})$. On the
other hand, $V(\mathcal{T})\subseteq V(\mathcal{P})$ is obvious. So,
$V(\mathcal{P})= V(\mathcal{T})$. Also, in terms of (a), $T$ is
connected. And, in terms of (c), $T$ has no cycles. Hence,
(\romannumeral2) holds.

In terms of (\romannumeral2) and (c), (\romannumeral3) holds.

Note that the path system $\mathcal{T}$ is the complete path system
on $[T,s]$. In terms of (\romannumeral3) and (a), we can easily know
that (\romannumeral4) holds.

The proof completes. $\hfill\square$

\section{Applications}
\label{sec:5}

This section shows the application of algorithm EDA and algorithm
EMBFA by providing few instances.

~\\
\textbf{Example 1.} Given a connected network with nonnegative
weight $(G,w)$ and source point $s\in V$. Let $\mathcal{P}$ be the
complete path system on $[G,s]$. Define
$d(P)=\sum\limits_{i=1}^{k}w((v_{i-1},v_{i})), \forall
P=(v_{0},v_{1},\cdots,v_{k-1},v_{k})\in \mathcal{P}, v_{0}=s$, which
is called the total weight path function of networks.
 It is obvious that $d$ is a path functional on $\mathcal{P}$. The problem
to find a path $P^\ast\in\mathcal{P}(v)$ such that
$d(P^\ast)=m_{d}(v), \forall v\in [V\setminus\{s\}]$, is called as
the classical single-source shortest path problem with nonnegative
weight (CSSSP-NW). See e.g. 7.1 of the monograph [11] (Korte and
Vygen 2000). Note that $w$ is nonnegative. We can easily know that
$d$ is nondecreasing and OP.
 So, by Theorem 2 and Theorem 3, the problem CSSSP-NW can
 be effectively solved by the
algorithm EDA and algorithm EMBFA, respectively.

~\\
\textbf{Remark 6.} From the example 1, we can know the following
facts clearly. (\romannumeral1) EDA and EMBFA respectively reduces
to DA and MBFA in the situation of the example. (\romannumeral2) Let
$P,P'\in\mathcal{P}$ and
$S(P)=P+(u,v),S(P')=P'+(u,v)\in\mathcal{P}$. Then
$d(S(P))-d(P)=d(S(P'))-d(P')=w((u,v))$ for the function $d$.
However, $f(S(P))-f(P)=f(S(P'))-f(P')$ may  not hold for a general
path function $f$. That is, we may not find an edge weight of graph
$G$ such that $f(P)=\sum\limits_{i=1}^{k}w((v_{i-1},v_{i})), \forall
P=(v_{0},v_{1},\cdots,v_{k-1},v_{k})\in \mathcal{P}$ for a general
path functional $f$. The two facts (\romannumeral1) and
(\romannumeral2) fully illustrate that EDA and  EMBFA respectively
extended DA and MBFA.

~\\
\textbf{Example 2.} For the example 1, change the nonnegative weight
 $w$ as a conservative weight, see e.g. Definition
7.1 of the monograph [11] (Korte and Vygen 2000). Then the problem
to find a path $P^\ast\in\mathcal{P}(v)$ such that
$d(P^\ast)=m_{d}(v), \forall v\in [V\setminus\{s\}]$, is called as
the classical single-source shortest path problem with conservative
weight (CSSSP-CW). Clearly, we can  easily know that $d$ is
conservative and OP. So, by  Theorem 3, the problem CSSSP-CW can be
effectively solved by Algorithm EMBFA.

~\\
\textbf{Example 3.} For the example 1, define also
$d(u,v)=\min\{d(P)|s(P)=u,t(P)=v,P\in\mathcal{P}_{nc}\}, \forall
u,v\in V$, which is called as the distance from $u$ to $v$ on graph
$G$. For given $(u',v')\in R(\mathcal{P}_{nc})$, define
$d_{G\backslash(u',v')}(u,v)$ as the distance from $u$ to $v$ on
graph $(V,[R\setminus\{(u',v')\}]), \forall u,v\in V$, which is
called as the detour distance from $u$ to $v$ on the case that the
road $(u',v')$ is blocked. Here $G\backslash(u',v')$ denotes the
graph $(V,[R\setminus\{(u',v')\}])$. In addition, $\forall(u',v')\in
R(\mathcal{P}_{nc})$, assume that $G\backslash(u',v')$ is connected,
namely $d_{G\backslash(u',v')}(u,v)<+\infty$ for each $(u,v)\in
[R\setminus\{(u',v')\}]$.

Define
$$\begin{array}{lcl}r(P)&=&0, P=(s,s);\\r(P)&=&\max\{d(P),
d_{G\backslash(v_{i-1},v_{i})}(s,v_{i})+d(P_{i}),d_{G\backslash(v_{k-1},v_{k})}(s,v_{k})|\\
 & &P_{i}=(v_{i},\cdots,v_{k-1},v_{k}),1\leq i\leq k-1\},\\
 & &\forall P=(v_{0},v_{1},\cdots,v_{k-1},v_{k})\in\mathcal{P}_{nc}, k\geq 1;\\
 r(P)&=&\max\{d(P)|P\in\mathcal{P}_{nc}\}+1,\forall P\in[\mathcal{P}\setminus\mathcal{P}_{nc}].
\end{array}$$
 We call
$r(P)$ as the risk of $P$. It is obvious that $r$ is a path function
on $\mathcal{P}$. The problem to find a path $P\in\mathcal{P}(v)$
such that $r(P)=m_{r}(v)$ for any $v\in
[V(\mathcal{P})\setminus\{s\}]$ is called the anti-risk path (ARP)
problem, the purpose of which is to find a path such that it has
minimum risk. See the literature [22] (Xiao et al. 2009).

Let $P=(v_{0},v_{1},\cdots,v_{k-1},v_{k})\in\mathcal{P}_{nc}$
 and
$S(P)=(v_{0},v_{1},\cdots,v_{k-1},v_{k},v_{k+1})\\=P+(v_{k},v_{k+1})\in\mathcal{P}_{nc},
k\geq 1$. Then we have
$$\begin{array}{lcl}& &r(S(P))\\&=&
\max\{d(S(P)),
d_{G\backslash(v_{i-1},v_{i})}(s,v_{i})+d(P_{i}),d_{G\backslash(v_{k},v_{k+1})}(s,v_{k+1})|\\
 & &\quad\quad\quad P_{i}=(v_{i},\cdots,v_{k},v_{k+1}),1\leq i\leq k\}\\
 &=&\max\{d(P)+w((v_{k},v_{k+1})),
d_{G\backslash(v_{i-1},v_{i})}(s,v_{i})+d(P_{i})\\& &\quad\quad\quad
+w((v_{k},v_{k+1})),d_{G\backslash(v_{k-1},v_{k})}(s,v_{k})+w((v_{k},v_{k+1})),\\&
& \quad\quad\quad
d_{G\backslash(v_{k},v_{k+1})}(s,v_{k+1})|P_{i}=(v_{i},\cdots,v_{k-1},v_{k}),
 1\leq i\leq k-1\}
\\&=&\max\{d_{G\backslash(v_{k},v_{k+1})}(s,v_{k+1}),w((v_{k},v_{k+1}))+\max\{d(P),\\& &\quad\quad\quad
 d(P_{i})+
 d_{G\backslash(v_{i-1},v_{i})}(s,v_{i}),d_{G\backslash(v_{k-1},v_{k})}(s,v_{k})\\& &\quad\quad\quad|P_{i}=(v_{i},\cdots,v_{k-1},v_{k}),1\leq i\leq k-1\}\}\\
&=&\max\{d_{G\backslash(v_{k},v_{k+1})}(s,v_{k+1}),w((v_{k},v_{k+1}))+r(P)\}\geq
r(P).\end{array}\eqno(8)$$ This shows $r$ is nondecreasing on
$\mathcal{P}_{nc}$.

Let also $P'\in\mathcal{P}_{nc}(v_{k}),
S(P')=P'+(v_{k},v_{k+1})\in\mathcal{P}_{nc}(v_{k+1})$. Assume
$r(P)\geq r(P')$. Then, from (8), we have
$$\begin{array}{lcl}r(S(P))&=&\max\{d_{G\backslash(v_{k},v_{k+1})}(s,v_{k+1}),w((v_{k},v_{k+1}))
+r(P)\}\\&\geq&\max\{d_{G\backslash(v_{k},v_{k+1})}(s,v_{k+1}),w((v_{k},v_{k+1}))
+r(P')\}=r(S(P')).\end{array}$$ This shows $r$ is SOP on
$\mathcal{P}_{nc}$.

Further, basing on the above two conclusions, we can easily show
that $r$ is nondecreasing and SOP on $\mathcal{P}$. So, by Theorem
2, the ARP problem can be effectively solved by algorithm EDA.

~\\
\textbf{Remark 7.} (\romannumeral1) Xiao et al. [22] (2009)
introduce the definition of the risk of a path, and the anti-risk
path (ARP) problem,  to finding a path such that it has minimum
risk. Suppose that at most one edge may be blocked, they also show
that the ARP problem can be solved in $O(mn+n^{2}log n)$ time.
Mahadeokar and Saxena [13] (2014) propose a faster algorithm to
solve the  ARP problem, by which the ARP problem can be solved in
$O(n^{2})$ time. (\romannumeral2) In example 3, for some technical
reason, the risk is defined in a slightly different manner from that
of Xiao et al. [22] (2009). Due to the cause of symmetry, in order
to conveniently understand the example 3 and the next example 4, the
path $P=(s,v_{1},\cdots,v_{k-1},v_{k},v)$ can be interpreted as the
path from $v$ to $s$. The ARP problem of example 3, which is
essentially similar to that of  Xiao et al. [22] (2009), can be
effectively solved by algorithm EDA. (\romannumeral3) For $r$ is not
proved to be OR, the ARP problem of example 3 may not necessarily be
solved by EMBFA. However, due that $r$ is SOP, it can be solved by
EDA. This fact fully demonstrates the advantages of SOP and EDA.

~\\
\textbf{Example 4.} For example 1, assume
$$d_{G\backslash(u',v')}(u,v)<+\infty, \forall(u',v')\in
R(\mathcal{P}_{nc}), \forall u,v\in V,$$ where
$d_{G\backslash(u',v')}(u,v)$ is the detour distance in example 3;
assume also $p\in(0,1)$.

Define first $c((s,s))=0$. Then, $\forall
P=(v_{0},v_{1},\cdots,v_{k})\in\mathcal{P}_{nc}, k\geq 1$, provided
that $c(F(P))$ has been defined, define
$c(P)=pd_{G\backslash(v_{k-1},v_{k})}(v_{0},v_{k})+w((v_{k-1},v_{k}))+c(F(P))$;
$\forall P\in[\mathcal{P}\setminus\mathcal{P}_{nc}]$, define
$c(P)=\max\{d(P)|P\in\mathcal{P}_{nc}\}+1$. Then $c$ is a path
functional on $\mathcal{P}$.

Let $P=(v_{0},v_{1},\cdots,v_{k-1},v_{k})\in \mathcal{P}_{nc}$, and
$S(P)=(v_{0},\cdots,v_{k}, v_{k+1})=P+(v_{k}, v_{k+1})\in
\mathcal{P}_{nc}$. Then we have
$$\begin{array}{lcl}c(S(P))&=&pd_{G\backslash(v_{k},v_{k+1})}(v_{k},v_{k+1})
+w((v_{k},v_{k+1}))+c(P)\\
&\geq& c(P).\end{array}$$ This shows $c$ is nondecreasing on
$\mathcal{P}_{nc}$.

Let also $P'\in\mathcal{P}_{nc}(v_{k}),
S(P')=P'+(v_{k},v_{k+1})\in\mathcal{P}_{nc}(v_{k+1})$. Assume
$c(P)\geq c(P')$.  Then,  we have
$$\begin{array}{lcl}
c(S(P))&=&pd_{G\backslash(v_{k},v_{k+1})}(v_{k},v_{k+1})+w((v_{k},v_{k+1}))+c(P)\\
&\geq&pd_{G\backslash(v_{k},v_{k+1})}(v_{k},v_{k+1})+w((v_{k},v_{k+1}))+c(P')\\
&=&c(S(P')).\end{array}$$ This shows $c$ is SOP on
$\mathcal{P}_{nc}$.

Further, we can easily show that $c$ is  nondecreasing and SOP. So,
by Theorem 2, the problem GSSSP with path function $c$ can be
effectively solved by algorithm EDA.

~\\
\textbf{Remark 8.} (\romannumeral1) The path function $c$ can be
interpret as follows. Suppose that at most one edge may be blocked.
$\forall P\in\mathcal{P}_{nc}$,
$c(P)=pd_{G\backslash(v_{k-1},v_{k})}(v_{k-1},v_{k})+w((v_{k-1},v_{k}))+c(F(P))$
denotes the cost that one goes to the point $s$ from the point
$v_{k}$ by train (or ship, or plane), among which, the term
$w((v_{k-1},v_{k}))+c(F(P))$ is the normal cost, while the term
$pd_{G\backslash(v_{k-1},v_{k})}(v_{k-1},v_{k})$ is the additional
cost, which is paid out due that one needs to change route when the
road $(v_{k-1},v_{k})$ is blocked. To some extent, $p$ is  the
probability that the road $(v_{k-1},v_{k})$ may be blocked.
(\romannumeral2) It is shown by example 3 and example 4 that path
functions is very useful in practice.

\section{Concluding remarks}
\label{sec:6} Due to the quick and extensive development in the
research  field of graph theories, and due to the quick and
extensive development in the research fields of networks and
functionals, especially with the emergence of the Anti-risk Path
Problem and the related research, the need to extend the total
weight path function of networks and the shortest path problems is
gradually increasing. Motivated by this trend, in the present
article, we have mainly done the following three aspects of work.

1. We proposed the definition of the path functional and several
definitions regarding
 the characteristics of the defined path
functional, especially the characteristic in orders, such as
increasing, order-preserving, so on; see Definition 1-Definition 3.
Further, based on the proposed definitions, we made some related
discussions on the properties of the path functional, see
Proposition 1-Proposition 6.

2. We introduced a kind of general single-source shortest path
problem (GSSSP) and designed two algorithms EDA and EMBFA to solve
it under certain conditions. Further, based on the discussions on
the properties of the path functional, we studied respectively the
attributes of EDA and EMBFA; see Theorem 2 and Theorem 3.

3. We further explained the significance of the defined path
functional and the two designed algorithms by several examples.

What we have done not only extends the Dijkstra algorithm  and the
Moore-Bellman-Ford algorithm, but also more profoundly reveals their
mechanism, which will greatly promote our understanding and applying
of the two algorithms. The discussions on the properties of the path
functional, not only support our designing and studying the two
extended algorithms, but also more profoundly reveals that the
contents of path functionals, in particular the content about the
order, are quite rich, depth and interest, which further shows that
there are many other research problems about path functionals.

It is an interesting topic for further research in the future to
explore other efficient algorithms and the applications of the
problem GSSSP. Moreover, to study other questions about path
functionals is also an interesting topic for further research in the
future. Finally, cordially hope that the present work can improve
the development of the researches and applications of shortest path
problems as well as other problems of combinatorial optimization.

\vskip .5cm
   \noindent{{\bf Acknowledgements}
The authors cordially thank the editor and the anonymous referees
    for their valuable comments and suggestions which lead to the improvement of
    this paper.


\newpage
\begin{thebibliography}{}

\bibitem{biblabel[AbuSalim]} AbuSalim, S. W. G., Ibrahim, R., Saringat, M. Z.,  et al.:
Comparative analysis between Dijkstra and Bellman-Ford algorithms in
shortest path optimization. In: IOP Conf. Series: Materials Science
and Engineering 917 (2020) 012077, 1-11 (2020)
https://doi.org/10.1088/1757-899X/917/1/012077

\bibitem{biblabel[Bellman]} Bellman, R. E.: On a routing problem. Quarterly of Applied
Mathematics, 16, 87-90 (1958)

\bibitem{biblabel[Dijkstra]} Dijkstra, E. W.: A note on two problems in connexion with graphs.
Numer Math, 1, 269-271 (1959)

\bibitem{biblabel[Deen]} Deen, M. R. Z. E., Aboamer, W.A., El-Sherbiny, H. M.: Explicit
Formulas for the Complexity of Networks Produced by New Duplicating
Corona and Cartesian Product. Journal of Mathematics, 2024, 9131329
(2024) https://doi.org/10.1155/2024/9131329

\bibitem{biblabel[Dinitz]} Dinitz, Y., Itzhak, R.: Hybrid Bellman-Ford-Dijkstra algorithm.
Journal of Discrete Algorithms, 42, 35-44 (2017)
http://dx.doi.org/10.1016/j.jda.2017.01.001

\bibitem{biblabel[Du]} Du, D. Z., Graham, R. L., Pardalos, P. M., et al.: Analysis of
greedy approximations with nonsubmodular potential functions. In:
Proceedings of the Nineteenth Annual ACM-SIAM Symposium on Discrete
Algorithms, 167-175 (2008)

\bibitem{biblabel[Feng]} Feng, G.: Finding k shortest simple paths in directed graphs: a
node classification algorithm. NETWORKS, 6-17 (2014)
https://doi.org/10.1002/net.21552

\bibitem{biblabel[Ford]} Ford, L. R.: Network flow theory. Paper P-923, The Rand
Corporation, Santa Monica, (1956)

\bibitem{biblabel[Hershberger1]} Hershberger, J., Suri, S.: Vickrey prices and shortest paths:
what is an edge worth?. In: Proceedings of the 42nd annual IEEE
symposium on foundations of computer science, 252-259 (2001)

\bibitem{biblabel[Hershberger2]} Hershberger, J., Suri, S., Bhosle, A.: On the difficulty of some
shortest path problems. In: Proceedings of the 20th symposium on
theoretical aspects of computer science, 343-354 (2003)

\bibitem{biblabel[Korte]} Korte, B., Vygen, J.: Combinatorial optimization theory and
algorithms.  Springer-Verlag, Berlin, (2000)

\bibitem{biblabel[Kov]} Kov\'{a}cs, L.: Classiffcation Improvement with Integration of
Radial Basis Function and Multilayer Perceptron Network
Architectures. Mathematics, 2025(13), 1471 (2025) https://
doi.org/10.3390/math13091471

\bibitem{biblabel[Mahadeokar]} Mahadeokar, J., Saxena, S.: Faster algorithm to find anti-risk
path between two nodes of an undirected graph. J Comb Optim, 27,
798-807 (2014)  https://doi.org/10.1007/s10878-012-9553-0

\bibitem{biblabel[Meng]} Meng, S. C. Y., Adnan, N., Sukri, S. S., et al.: An experiment
on the performance of shortest path algorithm. In: Knowledge
Management International Conference (KMICe) 2016, Chiang Mai,
Thailand, 7-12 (2016)

\bibitem{biblabel[Moore]} Moore, E. F.: The shortest path through a maze. Proccedings of
the international Symposium on the Theory of Switching, Part II.
Harvard University Press, Boston, 285-292 (1959)

\bibitem{biblabel[Murota]} Murota, K., Shioura, A.: Dijkstra's algorithm and L-concave
function maximization. Math. Program. Ser. A, 145, 163-177 (2014)
https://doi.org/10.1007/s10107-013-0643-2

\bibitem{biblabel[Nip]} Nip, K., Wang, Z., Nobibon, F. T., et al.:
 A combination of flow shop scheduling and the shortest
path problem. J Comb Optim, 29, 36-52 (2015)
https://doi.org/10.1007/s10878-013-9670-4

\bibitem{biblabel[Srivastava]} Srivastava, K.,  Tyagi, R.: Shortest path algorithm for
satellite network. www.ijird.com, 5(2), 438-445 (2013)

\bibitem{biblabel[Wu]} Wu, H.C.: Bernstein Approximations for Fuzzy-Valued Functions.
Mathematics, 2025(13), 2424 (2025) https://doi.org/10.3390/
math13152424

\bibitem{biblabel[Xu]} Xu, D.K., Tian, R.Q., Lu, Y.: Bayesian Adaptive Lasso for the
Partial Functional Linear Spatial Autoregressive Model. Journal of
Mathematics, 2022, 1616068 (2022)
https://doi.org/10.1155/2022/1616068

\bibitem{biblabel[Xia]} Xia, D. W., Shen, B. Q., Zheng, Y. L., et al.:
Abidirectional-a-star-based ant colony optimization algorithm for
big-data-driven taxi route recommendation. Multimedia Tools and
Applications, 83, 16313-16335 (2024)
https://doi.org/10.1007/s11042-023-15498-4

\bibitem{biblabel[Xiao]} Xiao, P., Xu, Y., Su, B.: Finding an anti-risk path between two
nodes in undirected graphs. J Comb Optim, 17, 235-246 (2009)
https://doi.org/10.1007/s10878-007-9110-4

\bibitem{biblabel[Zhang]} Zhang, H. L., Xu, Y. F., Wen, X. G.: Optimal shortest path set
problem in undirected graphs. J Comb Optim, 29, 511-530 (2015)
https://doi.org/10.1007/s10878-014-9766-5


\end{thebibliography}
\end{document}